\documentstyle [11pt,amsfonts]{article}

\newcommand{\ZZ}{{\Bbb Z}}
\newcommand{\NN}{{\Bbb N}}
\newcommand{\RR}{{\Bbb R}}
\newcommand{\QQ}{{\Bbb Q}}
\newcommand{\CC}{{\Bbb C}}
\newcommand{\CX}{{\Bbb C}^\times}

\newcommand{\BBox}{\rule{1.3ex}{1.3ex}}

\newcommand{\LB}{L_\Box}
\newcommand{\XB}{X_\Box}
\newcommand{\Sig}{\Sigma}
\newcommand{\sg}{\sigma}
\newcommand{\Ssig}{S_{\sigma}}
\newcommand{\usig}{U_{\sigma}}
\newcommand{\XNS}{X_{N,\Sigma}}
\newcommand{\XS}{X_{\Sigma}}
\newcommand{\NR}{N_{\RR}}
\newcommand{\bolde}{{\bf e}}
\newcommand{\bt}{{\frak t}}
\def\la{{\lambda_\alpha}}
\newcommand{\minus}{-\kern -0.2em}
\newcommand{\eqdef}{\stackrel{{\rm def}}{=}}

\newcommand{\vect}[3]{{({#1}_{#2},\dots,{#1}_{#3})}}

\newcommand{\ntup}[2]{#1,\dots,#2}
\newcommand{\gen}[2]{\langle#1,\dots,#2\rangle}

\newcommand{\ip}[2]{\langle#1,#2\rangle}
\newcommand{\inv}{{}^{-1}}
\newcommand{\trans}{{}^t\kern -0.2em}
\newcommand{\dual}{{}^{\vee}}
\newcommand{\sbot}{{}^{\bot}}
\newcommand{\sast}{{}^{\ast}}
\newcommand{\dast}{{}_{\ast}}

\newcommand{\trace}{{\rm trace}\,}
\newcommand{\ext}{{\rm ext}\,}
\newcommand{\orb}{{\rm orb}\,}

\newtheorem{dfn}{Definition}
\newtheorem{thm}{Theorem}
\newtheorem{prop}{Proposition}
\newtheorem{cor}{Corollary}
\newtheorem{lemma}{Lemma}

\newenvironment{proof}{{\em Proof:\/}}{$\BBox$\medbreak\par}

\begin{document}
\title{Lefschetz Fixed-Point Theorem and Lattice Points in Convex
Polytopes\thanks{M.S.Classification (1991): Primary: 52B20; Secondary: 14M25,
11H06, 11P21.} %alg-geom/9302003
}
 \author{Sacha Sardo Infirri \thanks{email:sacha@maths.oxford.ac.uk} \\
\\
Mathematical Institute\\
24-29 St. Giles'\\
 Oxford}
\date{\ }
\maketitle
\vskip -1cm
\begin{abstract}
      A simple convex lattice polytope
$\Box$ defines a torus-equivariant line bundle $\LB$ over a toric
variety $\XB.$

Atiyah and Bott's Lefschetz fixed-point theorem is applied to the torus action
on the $d''$-complex of $\LB$ and information is obtained about the lattice
points of $\Box$.

In particular an explicit formula is derived, computing the number of lattice
points and the volume of $\Box$ in terms of geometric data at its extreme
points. The results of Brion \cite{brion} are recovered and an elementary
convex geometric interpretation is given by performing Laurent expansions
similar to those of Ishida \cite{ishida}.
\end{abstract}

\setcounter{section}{-1}
\section{Introduction}
\subsection{The problem}
Let $\Box$ be a convex polytope all of whose vertices belong to a lattice $M$.
The question of calculating the number of points of
$M$ contained in $\Box$ is a well-known
one in convex geometry.  The oldest formula appears to be
Pick's classical result \cite{pick}, valid for arbitrary
polygons in 2 dimensions:
$$\#(\Box\cap M)={\rm Area\,}(\Box)+{1\over2}\#({\rm
boundary}(\Box)\cap
M)+1.$$
 Following Ehrhart's work on Hilbert polynomials, Macdonald
\cite{mac:lattice,mac:poly} subsequently generalised Pick's
formula to arbitrary $n$. His formula expresses
the volume of $\Box$ in terms of the number of lattice points
of its' multiples $k\Box$ for
finitely many integers $k$. Although these formulae are valid
for arbitrary (non-convex) polygons, they do not give any
convenient way of calculating either the volume, or the
number of
lattice points of $\Box$.
A review of this and other problems concerning lattice
points can be found
in \cite{hammer,erdos}.

{}From an elementary point of view, for large polytopes one
expects the volume
to be a good approximation to the
number of lattice points, so that one can imagine a general
formula of the
form
\begin{equation}\label{eq:RR}
\mbox{number of points = volume + correction terms}
\end{equation}
where the corrections terms are negligible in the large
limit. The formula we present here is however quite different in
nature.

\subsection{The results}

 Given a parameter $\zeta$, to each extreme point $\alpha$
of a simple convex polytope $\Box$, we associate a
rational number depending the local geometry of $\Box$ at $\alpha$.  Their sum
is independent of $\zeta$ and yields the number of lattice points in $\Box$
(Theorem \ref{thm:number}). Our main formula (Theorem \ref{thm:formula-sing})
is more general, since it expresses not just the number, but {\em which\/}
points of the lattice belong to the polytope, as a finite Laurent polynomial in
$n$ variables (the lattice points corresponding to the monomials via $m\mapsto
x^m=x_1^{m_1}x_2^{m_2}\cdots x_n^{m_n}$). I give an initial form of this using
the Lefschetz-fixed point theorem for orbifolds. By expanding in Laurent series
this is shown to be equivalent to another formulation (Theorem
\ref{thm:formula-sing-b}) given by Brion \cite{brion} which doesn't involve
cyclotomic sums. I use this form to calculate the number of lattice points. The
volume of $\Box$ is obtained by taking the leading order terms for finer and
finer subdivisions of the lattice. The Laurent series expansions extend
Ishida's \cite{ishida} and provide a convex geometric interpretation of the
formula (Theorem \ref{thm:chi-decomposition}). This in turn suggests a proof of
the formula involving {\em no toric geometry\/} --- only convex geometry and
elementary Laurent expansions. This could be considered as a variation of
Ishida's proof \cite{ishida} based on the contractibility of convex sets.

This paper is an amplification of my 1990 transfer dissertation at Oxford
university \cite{sacha}. This was originally written whilst I was unaware of
Michel Brion's 1988 paper \cite{brion}, where a toric approach is used to
calculate the number of lattice points. There has  also been a paper by Ishida
\cite{ishida} where similar Laurent expansions similar to mine are performed.
This is the revised version of my original which takes these works into
account. Let me briefly mention their relationship to this paper.

Brion relies on the Lefschetz-Rieman-Roch theorem for equivariant K-theory
\cite{BFQ} and obtains theorem \ref{thm:formula-sing-b}. He calculates the
number of lattice points by subdividing the tangent cones into basic cones.
The formula that I obtain using the Lefschetz fixed point theorem involves
instead cyclotomic sums for the action of the finite quotient group. By
extending Ishida's Laurent series expansions \cite{ishida} in section
\ref{sec:laurentexpansions} of this article, I prove that the two are
equivalent, and provide a combinatorial interpretation of the formula. It is
also not necessary for me to subdivide the tangent cones in order to obtain a
formula for the number of lattice points.

\subsection{The method}

Our main tool is the theory of toric varieties. This
associates a holomorphic line bundle $\LB$ over a complex
orbifold $\XB$ to any n-dimensional simple  polytope $\Box$
on a lattice $M$. The variety comes equipped with the action of an algebraic
$n$-torus $T_N$ (the character group of $M$) and $\LB$ is equivariant with
respect to this action. Its cohomology  is trivial in positive dimension,
whereas its space of sections is naturally isomorphic to a vector space
generated by the lattice points in $\Box$.

In \cite{bern,khov}, the Rieman-Roch theorem is used to
calculate the number of lattice points in $\Box$. This yields a
formula similar to equation (\ref{eq:RR}) above. The
problem with this approach, however, is that the correction
terms are not readily computable.

In this paper I follow an idea of Atiyah and exploit the
torus action. I apply Atiyah \& Bott's Lefschetz fixed point theorem
\cite{ab:lefI} --- suitably extended to orbifolds \cite{kawasaki} --- to the
(geometric endomorphism induced by the) action of $t\in T_N$ on the
$d''$-complex of $\LB$. The $d''$-complex is elliptic \cite{ab:lefII} and its
cohomology groups are (canonically isomorphic to) those of $(\XB,\LB)$. The
fixed points of the torus action on $\XB$ correspond to the extreme points of
$\Box$. The Lefschetz theorem in this case expresses the equality between the
Lefschetz number (an element of $\CC[M]$)
and the sum of the indexes $\nu_\alpha$ for $\alpha$ in the set of extreme
points. The $\nu_\alpha$ define elements of $\CC(M)$. The formula I obtain
initially involves sum over the characters of the finite abelian groups which
charaterise the singularities at the points $P_\alpha\in\XB$ corresponding to
$\alpha\in\ext\Box$. By studying characteristic series for cones in section
\ref{sec:laurentexpansions} I eliminate the summation over group elements.

If one restricts $t$ to the one-parameter subgroup of the torus determined by
an element $\zeta$ of its' Lie algebra one obtains an equality between a
polynomial and a sum of rational functions in one variable. When $t\to 1$ the
polynomial tends to the the number of lattice points of $\Box$, and this is
given by the sum of the constant terms in the one variable Laurent series for
the rational functions: this gives theorem \ref{thm:number}. By identifying the
coefficient of the leading order terms in the asymptotic expansions of the
formula for  submultiples of the lattice --- the `classical limit' in quantum
terminology --- I derive a formula for the volume of $\Box$ in Theorem
\ref{thm:volume}.

I review the toric geometry results I shall need in the first part of this
paper. The reader who is familiar with the notation in Oda \cite{oda} can {\tt
GOTO PART II}, which contains the application proper.

\subsection{Acknowledgments}
I would like to thank Michael Atiyah and Peter Kronheimer for their stimulating
ideas and encouraging support. Thanks also to Frances Kirwan for her
suggestions and to Mark Lenssen,  Jorgen Andersen and Jorge Ramirez-Alfonsin
for interesting discussions. I was supported by a Rhodes Scholarship while I
did this research.

\subsection{Notation}
\label{subsec:notation}

Throughout this paper, let $N\cong \ZZ^n$ denote an n-dimensional
integral lattice, $M\cong\hom_\ZZ(N,\ZZ)$ it's dual and
$\NR=N\otimes_{\ZZ}\RR$ its' associated real vector space. The
complex torus $N\otimes_{\ZZ}
\CX\cong\hom_{\ZZ}(M,\CX)$ is denoted $T_N$ and the compact
real sub-torus $N\otimes S^1\subset N\otimes \CX$ is denoted $CT_N$.

If $A$ is any commutative ring with identity and $S$ any additive semi-group,
we write $A[S]$ for the group algebra of $S$ with   coefficents in $A$; this is
generated by elements $\bolde(s)$ for $s\in S$ satisfying the relations
$\bolde(s)\bolde(s')=\bolde(s+s')$. We write $A(S)$ for its total quotient ring
(i.e., its' field of fractions if $A=\CC$).

Occasionally I choose coordinates $t_i$ for $T_N$.
This is equivalent to choosing generators $n_i$ for $N$. I denote the dual
generators by $m^j\in M$. Then if $\alpha\in M$ and $z \in T_N$ have
coordinates $\vect\alpha1n$ and $\vect z1n$ with respect
to the appropriate bases we have
$$\alpha(z)=z^\alpha=z_1^{\alpha_1}z_2^{\alpha_2}\cdots
z_n^{\alpha_n}.$$
This identifies $\CC[M]$ with the Laurent polynomials in the variables $t_i$.
\newpage
\part{Toric Geometry}

The theory of toric varieties establishes correspondances
between convex geometry in $n$ real dimensions and  the
geometry of compactifications of $n$-dimensional complex
tori. I refer to \cite{kempf,oda,danilov}.

Briefly, there is a functor that associates, to a pair
$(N,\Sigma)$ (where $\Sigma$ is a fan in $N$),
an irreducible normal Hausdorff complex analytic space $\XNS$.
A convex polytope $\Box$ in $M$ determines
a unique fan $\Sig$ in $N$, and we set \(\XB = \XNS \).  The
polytope contains more information than simply its cone
structure, and this determines a piecewise linear function
$h=h_{\Box}$ on the support $|\Sig|\subset \NR$ of $\Sig$. This
corresponds under the functorial construction above to an
equivariant line bundle $L_h$ on $\XNS$, which we denote by
$L_{\Box}$.

\section{Cones and Affine Toric Varieties}

\subsection{Cones}

Let $V$ denote a vector space and $V^*$ its dual.

A {\em cone\/} in a vector space $V$ is a finite intersection
of half-spaces in $V$. Cones are always convex and polyhedral. I shall take
them to be also strongly convex, namely such that they do not contain any
proper subspace of $V$.

For $\ntup{v_1}{v_k} \in \NR$,
let $\gen{v_1}{v_k}$ denote the smallest cone containing
$\ntup{v_1}{v_k}$. Any cone is generated in this way.  A cone
is said
to be {\em simplicial\/} if it can be generated by linearly
independent elements of $\NR$. If it can be generated by part
of a $\ZZ$-basis of $N$, then the cone is called {\em
basic\/}.  Finally, a cone is
said to be {\em integral\/} with respect to $N$ if it can be
generated by elements of $N$.
When we speak of a {\em cone in a lattice\/} $N$ we mean a
cone in  $\NR$ which is integral with respect to $N$. I only
consider such cones henceforth.

 The {\em dimension\/} of a cone is the dimension of the subspace it
 generates. By the {\em interior\/} of a cone we usually mean the relative
interior in the subspace it generates.

\subsection{Duality}
Given a subset $A\subset V$  its {\em dual\/}
$A\dual\subset V\sast$ is defined by:
\[ A\dual=\{\theta\in V\sast : \forall v\in V,
\ip{\theta}{v} \ge 0\}.\]

\begin{prop}
\label{prop:duality}
 The dual of a cone (respectively, a simplicial cone, a basic cone,
 or an integral cone) is a cone (respectively a simplicial cone,
 a basic cone, or an integral cone). Moreover, for
any cone $\sg$ we consider, we have $\sg\dual\dual=\sg$.
\end{prop}
 For a proof of all the results regarding cones, see
\cite{rockaf}.
A summary of the results I require will be found in
\cite{oda}.

\subsection{Affine Toric Varieties}

Let $\sg$ be a cone in $N$. Recall  \cite[Prop.\ 1.1]{oda} that
the subset of $M$ given by
\[ \Ssig = M\cap \sg\dual \] is finitely generated as an
additive semigroup, generates $M$ as a group, and is
saturated. Such semigroups are in one-one correspondance
with cones in $N$.

Denote by $\usig=U_{N,\sigma}$ the set of semigroup
homomorphisms from
$(\Ssig, +)$ to $(\CC, \cdot)$, namely
\[\usig=\{u:\Ssig\to\CC: u(0)=1,
u(m+m^{\prime})=u(m)u(m^{\prime}),\forall m,m^{\prime}
\in \Ssig\}.\]
This can be given the structure of an n-dimensional
irreducible
normal complex analytic space by choosing generators
 $\ntup{m_1}{m_p}$ for $\Ssig$ and embedding $\usig$ in
$\CC^p$ via the
evaluation maps ${\bf ev}(m_i) :  u\mapsto u(m_i)$ on the
generators $m_i$.
The structure is inherited from the usual structure on
$\CC^p$ and is
independent of the generators chosen.

In other words, $\usig$ is just equal to the (set of points of the) affine
scheme $\mbox{Spec}(\CC[\Ssig])$. Identifying $\usig$ with its $\CC$-points
corresponds to identifying ${\bf ev}(m)$ with $\bolde(m)$. I spend little
effort making the distinction. The following proposition is easy to show
\cite[Th. 1.10]{oda}:

\begin{prop}
\label{prop:non-sing}
 The variety $\usig$ is non-singular if and only if $\sigma$
is basic.
\end{prop}

\section{Fans and General Toric Varieties}

\subsection{Faces, Fans and Gluing}

Let $\sg$ be a cone in $N$.

\begin{dfn}
A {\em face\/} of $\sg$ is a subset of the form
$\sg\cap\{m_0\}\sbot$,
where $m_0\in M=\hom(N,\ZZ)$ is non-negative on $\sg$. A
face of a cone is also a cone.
\end{dfn}

We immediately have:

\begin{lemma}
\label{lemma:open}
 If $\tau$ is a face of $\sigma$ then, for some $m_0\in M$,
we have
 \[ U_{\tau} = \{u\in\usig : u(m_0)\ne 0\},\]
so that $U_{\tau}$ is naturally an open subset of $\usig$.
\end{lemma}

Given this, one constructs collections of cones (called
{\em fans\/}) which have the property that their
corresponding varieties fit together in a natural way:

\begin{dfn}
A {\em fan\/} in $N$ is a collection $\Sigma=\{\sigma:
\sigma \mbox{ a cone
in }N\}$ satisfying the following conditions:
\begin{itemize}
\item if $\tau$ is a face of $\sg$ and $\sg \in\Sig$, then
$\tau\in \Sig$.
\item $\sg\cap\sg^{\prime}$ is a face of both $\sg$ and
$\sg^{\prime}$, for all
$\sg,\sg^{\prime}\in\Sig$.
\end{itemize}
\end{dfn}

The set of cones of $\Sig$ of
dimension $k$ is
called the {\em k-skeleton\/} of $\Sig$ and is denoted
$\Sig^{(k)}$.
The union of all the cones of $\Sig$ is called the {\em
support\/}
of $\Sig$ and is denoted $|\Sig|\subset \NR$.

\begin{thm}
\label{thm:general_toric}
The {\em toric variety\/} associated to $(N,\Sig)$ is the
space obtained by gluing together the affine varieties
$U_{N,\sg}$ for $\sg\in\Sig$, using lemma \ref{lemma:open}.
It is an n-dimensional Hausdorff complex analytic space
$\XNS$ which is irreducible and normal \cite[Theorem
1.4]{oda}. It is compact if and only if $\Sig$ is {\em
complete\/}, namely if and only if  $|\Sig|=\NR$.
\end{thm}

\subsection{The torus action}

The torus $T_N$ acts on $\usig$ by \((t\cdot u)(m)=t(m)u(m)\),
and this gives an action on $\XNS$. For $\sg=\{0\}$, one has
$U_{\{0\}}=T_N$, and the action coincides with group
multiplication on the torus.

The $T_N$-orbits on $\XNS$ are given by the
quotient algebraic tori
\begin{equation}
\label{eq:orb}
\mbox{orb}(\tau)=\hom_{Z}(M\cap\tau\sbot,\CX),
\end{equation}
for each $\tau\in\Sig$. The orbit corresponding to $\tau$
has
dimension equal to the codimension of $\tau$ in $\NR$. It is
also
easy to see that
$\usig$ decomposes as the disjoint union of the orbits
corresponding to its faces, and that $\mbox{orb}(\sg)$ is the
only
closed orbit in
$\usig$. I record a special case of this for later use:
\begin{lemma}
\label{lemma:fixpts}
The fixed points of the $T_N$ action on $\XNS$ are in one-
one
correspondance with the orbits
$\mbox{orb}(\sigma)\in\usig$, for the cones
$\sg$ in the $n$-skeleton $\Sig^{(n)}$.
\end{lemma}

\subsection{Functoriality}

Recall the following characterisation of toric varieties:
\begin{quote}
\em
 $X$ is a toric variety if and only if it is an irreducible
normal variety, locally of finite type over $\CC$, with a
densely embedded torus whose action on itself extends to
the whole variety.
\end{quote}

 The assignment $(N,\Sigma) \mapsto \XNS$ is a functor of
categories:

\begin{dfn} A {\em map of fans\/}
\(\phi:(N^{\prime},\Sig^{\prime})\to(N,\Sig)\) is
a $\ZZ$-linear homomorphism \(\phi:N^{\prime}\to N\)
whose scalar
extension \(\phi_{R}:N^{\prime}_{R}\to N_{R}\) satisfies the
following property: for each $\sg^{\prime}\in\Sig^{\prime}$,
there
exists
$\sg\in\Sig$ such that $\phi_{R}(\sg^{\prime})\subset\sg$.
\end{dfn}

\begin{thm} \cite[page 19]{oda} A map of fans
\(\phi:(N^{\prime},\Sig^{\prime})\to(N,\Sig)\) gives
rise to a holomorphic map
\[\phi\dast:X_{N^{\prime},\Sigma^{\prime}}\to\XNS\]
whose restriction to the open subset $T_{N^{\prime}}$
coincides
with the
homomorphism of algebraic tori
\(\phi_{\CX}:N^{\prime}\otimes_{\ZZ}\CX\to
N\otimes_{\ZZ}\CX.\)
Through
this homomorphism, $\phi\dast$ is $(T_{N^{\prime}}, T_N)$-
equivariant.
Conversely any holomorphic map $\psi:X'\to X$ between toric
varieties which restricts to a homomorphism $\chi: T'\to T$
on the algebraic tori $T'$ and $T$ in such a way that $\psi$
is $\chi$-equivariant corresponds to a unique
$\ZZ$-linear homomorphism
\(f:N^{\prime}\to N\) giving rise to a map of fans
\((N^{\prime},\Sig^{\prime})\to(N,\Sig)\)
such that $f\dast=\psi$.
\end{thm}

\subsection{Finite Quotients}
I will be interested in the case when $N^{\prime}$ is a
$\ZZ$-submodule of $N$ of finite index and
$\Sig^{\prime}=\Sig$. I
write $X^{\prime}$
and $X$ for the corresponding varieties:

\begin{prop} With the data as above, $X^{\prime}\to X$
coincides
\label{prop:quotient}
with the projection of $X^{\prime}$ with respect to natural
action of the
finite group
\[K= N/N^{\prime} \cong\hom_{Z}(M^{\prime}/M,\CX)=
\ker[T_{N^{\prime}}\to T_N].\]
\end{prop}
\begin{proof}
\cite[Cor. 1.16, p.22]{oda}
\end{proof}

\section{Toric Varieties, Equivariant Line Bundles and
Convex Polytopes}

\subsection{Polytopes}
\label{subsec:polytopes}
Recall first some basic notions of convex geometry.

A {\em convex polytope\/} $\Box$ in a vector space $V$ is a
bounded intersection of a finite number of affine half-
spaces of
$V$. The set of extreme points of $\Box$ is denoted
$\ext\Box$. Since $\Box$ is bounded, it is equal to the convex
hull of $\ext\Box$.

By a {\em polytope on the lattice\/} $M$ we mean a polytope
in
$M_\RR$ such that $\ext\Box\subset M$. Suppose $\Box$ is
such a
polytope, and let $\alpha$ be an extreme point. I define the
{\em (tangent) cone of $\Box$ at $\alpha$\/} to be the
cone $C_\alpha$ in $M$ given by:
\begin{equation}
\label{eq:calpha}
C_\alpha=\RR_{\ge 0}(\Box-\alpha)=\{r(v-
\alpha):r\ge0,v\in\Box\}.
\end{equation}

Let $\lambda_{\alpha}^{i}, i=1,\dots,k$ be the shortest
generators
for $C_{\alpha}$ which belong to the lattice $M$. I call
these the
{\em edges of $\Box$ emanating from}\/ $\alpha$, or simply the {\em
edge vectors for
$\Box$ at $\alpha$}. If $C_\alpha$ is simplicial (respectively,
 basic),  then $\Box$ is called {\em simple\/} (respectively,
{\em basic}) at $\alpha$. Henceforth, all the polytopes we
consider are convex, integral
and simple at all extreme points. They may be non-basic.

\subsection{Toric Varieties Defined by Polytopes}

\subsubsection{The Fan Defined by a Polytope}

The construction of $C_\alpha$ described in the previous
section can be generalised to show that a polytope $\Box$ in
$M$ defines a complete
fan in $N$. To each face $\Gamma$ of $\Box$ we associate the
cone
$C_{\Gamma}$ in $M$ defined by
\[ C_{\Gamma}=\RR_{\ge 0}(\Box-m_{\Gamma}),\]
where $m_{\Gamma}$ is any element of $M$ strictly in the interior of the
face $\Gamma$. If $F=\{\alpha\}$ we set $C_{\{\alpha\}}=C_\alpha$, as defined
previously in equation (\ref{eq:calpha}). Taking duals one obtains a collection
of cones in $N$
\[\Sig_\Box=\{\sigma_{\Gamma}=C\dual_{\Gamma} :
\Gamma
\mbox{ a face of }\Box\}.\]
One has the following easy lemma:

\begin{lemma}
\label{lemma:fanBX}
 $\Sig_{\Box}$ is equal to the fan consisting of the
cones
$\sigma_{\alpha}=C_{\alpha}\dual$, for $\alpha\in\ext\Box$
and
all their faces. It is complete, and its n-skeleton is
\(\Sig^{(n)}=\{\sigma_{\alpha}: \alpha\in\ext\Box\}\)
\end{lemma}

\subsubsection{The Variety Defined by a Polytope}

I define $\XB$ to be $\XS$, for $\Sig=\Sig_{\Box}$. By
\cite[Theorem
2.22]{oda}, $\XB$ is an orbifold (i.e., it has at worst quotient singularities)
if $\Box$ is simple.
\begin{prop}
\label{prop:XBaction}
The variety $\XB$ is compact, and is covered by affine
pieces
$$U_{\alpha}=U_{\sigma_{\alpha}}=\mbox{Spec}(\CC[M\cap
C_\alpha]),$$
for $\alpha\in\ext\Box,$ each containing a unique $T_N$--fixed
point $P_{\alpha}=\orb(\sigma_\alpha)$ (see equation
(\ref{eq:orb})).
Furthermore, when $U_{\alpha}$ is non-singular, the weights
of
the $T_N$ action on the tangent space
$T_{P_{\alpha}}U_{\alpha}$ are given by the edges vectors
for $\Box$ at $\alpha$.
\end{prop}
\begin{proof}
 The first claim follows directly from theorem
 \ref{thm:general_toric} and lemmas
\ref{lemma:fixpts}
and \ref{lemma:fanBX}. For the second part, observe (prop.\
\ref{prop:non-sing} and \ref{prop:duality}) that $U_{\alpha}$
is
non-singular if and only if the edge vectors at $\alpha$
generate
$M$ as a group. The semigroup $C_\alpha$ is then free on
these generators. They correspond to the weights of $T_N$ on
$U_\alpha$, and hence, by linearity, to the weights on
$T_{P_{\alpha}}U_{\alpha}.$
\end{proof}

\subsection{Equivariant Line Bundles}

The polytope $\Box$ contains more information than the fan
$\Sig_\Box$. This extra information turns out to be exactly
what
one needs to specify a $T_N$--equivariant  line bundle $\LB$
over
$\XB$.

\subsubsection{Line Bundles and Piecewise Linear Functions}
\label{subsub:LBPLF}
In general (equivalence classes of) equivariant line bundles
over
$\XNS$ are in one-one correspondence with the space
$PL(N,\Sig)$
of {\em piecewise linear functions} on $(N,\Sig)$, namely
functions \(h:|\Sig|\to\RR\) that are linear on each
$\sg\in\Sig$
and which take integer values on the integer points of
$|\Sig|$.

Defining an element $h\in PL(N,\Sig)$ involves, by
definition,
specifying an element $l_\sigma\in M$ for each
$\sigma\in\Sig$
such that $h(n)=\ip{l_\sigma}{n}$ for all $n\in\sigma$.
These
elements determine a line bundle $L_h$
equiped with a $T_N$-action and whose projection
$L_h\to\XS$ is equivariant
with respect to that action. Note that in general, the
elements $l_\sigma$ are not uniquely determined by $h$, but
different choices give rise to equivariantly equivalent
bundles.

 The bundle $L_h$ is defined to be trivial over the varieties
$U_\sigma$, with transition functions given by
\begin{equation}
g_{\tau\sigma}(x)=\bolde(l_\sigma-l_\tau)(x).
\label{eq:trans}
\end{equation}
The action of $T_N$ on the piece
$U_\sigma\times\CC\subset L_h$
is defined by
\begin{equation}
t(x,c)=(tx,\bolde(-l_\sigma)(t)c). \label{eq:actionL}
\end{equation}

\subsubsection{Cohomology}

The cohomology groups for equivariant line bundles
decompose
under the action of $T_N$ into weight spaces, and can be
expressed as a direct sum (see \cite[Th. 2.6]{oda}):
\[H^q(\XS,{\cal O}_{\XS}(L_h))=\oplus_{m\in M}
H^q_{Z(h,m)}(N_R,\CC)\bolde(m),\]
where $Z(h,m)=\{n\in N_R:\ip{m}{n} \ge h(n)\},$ and
$H^q_{Z(h,m)}(N_R,\CC)$ denotes the $q$-th cohomology group of
$N_R$
with support in $Z(h,m)$ and coefficients in $\CC$.

\paragraph{The Line Bundle $\LB$}

The polytope $\Box$ defines a piecewise linear function
$h_\Box$
on $\Sig_\Box$ by putting $l_{\sigma_{\alpha}}=\alpha$
(and
$l_\sigma=\alpha$ for the faces $\sigma$ of
$\sigma_{\alpha}$).
The corresponding bundle is denoted  $\LB$.
Its cohomology is given by \cite[Cor. 2.9]{oda}
\begin{equation}
\label{eq:coho}
H^q(\XB,{\cal O}_{\XB}(\LB))=\left\{
\begin{array}{ll}
\CC[M]_\Box = \oplus_{m\in M\cap\Box} \CC\bolde(m) & \mbox{if $q=0$} \\
0 & \mbox{otherwise}
\end{array}\right.
\end{equation}

\newpage
\part{The Polytope Formula}

\section{The Lefschetz Fixed-Point Theorem}

Recall \cite[Theorem 4.12]{ab:lefII} the following
application of
the Lefschetz fixed point theorem to the case of holomorphic
vector bundles:

\begin{thm}
\label{thm:lefschetz}
 Let $X$ be a compact complex manifold $X$, $F$ a
holomorphic
vector bundle over $X$, $f:X\to X$ a holomorphic map with
simple
fixed points and $\phi:f\sast F\to F$ a holomorphic bundle
homomorphism. Let $L(T)$ be the Lefschetz number of the
endomorphism $T$ of the $d''$-complex of $F$:
\[L(T)=\sum(-1)^q\mbox{trace}H^q T|_{H^q(X;F)}.\]
 Then
\(L(T)=\sum_{P=f(P)}\nu_{P}\), where
\[\nu_{P}={{\mbox{trace}_{\CC}
\phi_{P}}\over{\mbox{det}_{\CC}(1-df_{P})}}.
\]
\end{thm}

(Recall that since $P$ is a fixed point, $\phi_{P}$ and
$df_{P}$ are
endomorphisms of $F_{P}$ and $T_{P}X$ respectively.)

\subsection{Application}

       I apply this to the case where $X=\XB$, $L=\LB$ and
$f: X \to X$ is given by the action of a non-trivial element of $t\in T_N$. The
fixed points are simple and are given by
$P_{\alpha}=\mbox{orb}(\sigma_{\alpha})\in
U_{\sigma_{\alpha}}$, for $\alpha\in\ext\Box$.
The bundle homorphism $\phi_{t}: t\sast\L\to L$ is given by
the action of $-t$ (recall that $T_N$ acts on line bundles).
The cohomology groups are all zero, except $H^0(\XB,\LB)$ which is isomorphic
to the subspace $\CC[M]_\Box$ of $\CC[M]$ determined by $\Box$.

In this context, the Lefschetz number is an element of $\CC[M]$ and the indexes
$\nu_\alpha=\nu_{P_\alpha}$ are elements of $\CC(M).$ (As we shall see in
section~\ref{sec:laurentexpansions}, they are characteristic functions for the
tangent cones to $\Box$.)

\begin{lemma}
\label{lemma:Laction}
We have
\[\trace(\phi_{t})_{P_\alpha} = \alpha(t).\]
\end{lemma}
\begin{proof}
Recall (equation
 (\ref{eq:actionL})), that $t$ acts on the fibres of $L$ over
$\usig$ by
multipication by $\bolde(-l_{\sigma})(t)$, where
$l_{\sigma}$ are
the elements of $M$ corresponding to $L$ as in
\ref{subsub:LBPLF}. In the present case, at a fixed point
$P_\alpha\in U_{\alpha}$ we have
$l_{\sigma_{\alpha}}=\alpha$, so $\phi_t$ acts
by $\bolde(-\alpha)(-t)=\alpha(t)$.
\end{proof}

In the case of a basic polytope $\Box$ in $M$, applying
Theorem \ref{thm:lefschetz} directly one obtains:

\begin{thm}
\label{thm:non-sing}
For a basic simple convex polytope $\Box$ in $M$, we have
\begin{equation}\label{eq:sum-non-sing}
\sum_{m\in\Box}m(t)=\sum_{\alpha\in \ext\Box}
\nu_\alpha(t)
\end{equation}
where
\begin{equation}\label{eq:nu-non-sing}
\nu_\alpha(t)=\sum_{\alpha\in \ext\Box}
{\alpha(t)
\over
(1-{\lambda_{\alpha}^1}(t))
\cdots
(1-{\lambda_{\alpha}^n}(t))},
\end{equation}
the vectors
$\ntup{\lambda_{\alpha}^1}{\lambda_{\alpha}^n}$ are the
 edge vectors of $\Box$ at $\alpha$.
\end{thm}
\begin{proof}
The decomposition of $H^0(\XB;{\cal O}_{\XB}(\LB))$ given by
equation  (\ref{eq:coho}) shows that the left-hand side of
equation
 (\ref{eq:sum-non-sing}) is equal to the
Lefschetz number of the endomorphism induced by the
action of $t$.
Lemma \ref{lemma:Laction} and Proposition
\ref{prop:XBaction}
yield equation  (\ref{eq:nu-non-sing}).
\end{proof}

\subsection{The Lefschetz Fixed-Point Theorem for Orbifolds}

In \cite{kawasaki} the Lefschetz formula is generalised to
orbifolds (also known as V-manifolds), using zeta-function
techniques.
As I do not need the full power of this approach, I
present an alternative  more elementary argument.

The Lefschetz fixed-point formula is essentially local in
nature, the formula for the multiplicities $\nu_\alpha$ only
involving
the properties of $f$ and $\phi$ at the point $P_\alpha$.
This fact is
clearly apparent in Atiyah and Bott's proof in \cite{ab:lefI}
(see their
remarks at the beginning of section 5, and Proposition 5.3).
To extend the
formula to orbifolds, it is sufficient therefore to extend it
to global
quotient spaces, of the form $X=X'/K$.

\begin{prop}
\label{prop:Lef-Quotient}
Suppose   that a finite abelian group $K$ acts on a smooth
$X'$ and  equivariantly on a
holomorphic bundle $F'$ over $X'$. Let $f':X'\to X'$ and
$\phi':f'\sast F'\to F'$ be as in Theorem \ref{thm:lefschetz},
and suppose
they are $K$-equivariant. Denote by $L'(T')$ the Lefschetz
number of the
corresponding endomorphism $T'$ of the $d''$-complex of
$F'$. Because of the
$K$-equivariance, we can
define $X=X'/K$, $f:X\to X$, $F=(F')^K$, $\phi:f\sast F\to F$
and the
corresponding Lefschetz number
\[L(T)=\sum(-1)^q\trace H^qT|_{H^q(X;F)}.\]
Then we have
\begin{equation}\label{eq:L-Quotient}
L(T)={1\over{|K|}}\sum_{k\in K}L'(k\circ T).\end{equation}
\end{prop}
\begin{proof}
Note that since $T$ determines an endomorphism of the primed complex, it
makes sense to write
$L'(T)$. The claim then follows by applying the following
easy lemma of linear
algebra, recalling that $H^q(X;F)$ is just the $K$-invariant
part of
$H^q(X';F')$.
\end{proof}

\begin{lemma}
Suppose we have a linear action of a
finite abelian group $K$ on a finite dimensional vectorspace
$V$, commuting with an endomorphism $T$ of $V$. Denote by
$V^K$
the $K$-invariant subspace of $V$. Then $T$ is an
endomorphism
of $V^K$ and we have
$$\trace T|_{V^K}={1\over{|K|}}\sum_{k\in K}
\trace (k\circ T)|_V .$$
\end{lemma}
\begin{proof}
Define $P$ to be the following endomorphism of $V$:
$$Pv = {1\over{|K|}}\sum_{k\in K} k\cdot v.$$
Then $P^2=P$, so $P$ is the projection $V\to V^K$. Since $T$
commutes with $P$, it follows that $T$ respects the
decomposition
$V=V^K\oplus \ker P.$
 Furthermore we have $$\trace T|_{V^K}=\trace TP|_V  =
\trace
PT|_V ,$$
 so the result follows.
\end{proof}

Now, given a general orbifold $X$, at each point $P\in X$,
choose a {\em local model\/}  $(U'_P,f'_P,K_P,L'_P)$ as follows:

Let $U_P$ be an $f$-invariant neighbourhood of $P$ in $X$ and
$U'_P$ be a smooth cover with an action of a finite group $K_P$,
free away from $P$, such that $U_P=U'_P/K_P$. Thus $X$ has a quotient
singularity of type $K_P$ at $P$. Let $f'_P:U'_P\to U'_P$ be a
$K_P$-equivariant lifting of $f|_{U_P}$. A line bundle $L$ over $X$ is
understood to be an invertible sheaf $L$  over $X$ such that for any $P\in X$
with local model $(U'_P,f'_P,K_P)$, there exists a line bundle $L'_P\to U'_P$
such that $L|_{U_P}=(L'_P)^{K_P}$.

With these definitions, our remarks at the beginning of the
section and Proposition
\ref{prop:Lef-Quotient} imply the following:

\begin{thm}
\label{thm:lefschetz-orbifold}
 Let $X$ be a compact complex orbifold $X$, $F$ a
holomorphic
vector bundle over $X$, $f:X\to X$ a holomorphic map with
simple
fixed points and $\phi:f\sast F\to F$ a holomorphic bundle
homomorphism. Let $L(T)$ be the Lefschetz number of the
endomorphism $T$ of the $d''$-complex of $F$:
\[L(T)=\sum(-1)^q\mbox{trace}H^q T|_{H^q(X;F)}.\]
 Then
\(L(T)=\sum_{P=f(P)}\nu_{P}\), where
\[\nu_{P}={1\over |K_P|} \sum_{k\in K_P} {{\mbox{trace}_{\CC}
(k\circ\phi'_{P})}\over{\mbox{det}_{\CC}(1-(k\circ df')_{P})}},
\]
and $\phi', f'$ are lifts for $\phi, f$ respectively, in the same
spirit as that of the local models above.
\end{thm}

\subsection{Singular Case}

Suppose that $\Box$ is not basic relative to $M$ at
$\alpha$. Then $X=\XB$ has a singularity at the point $P_\alpha$. Let
$C_\alpha$ be the cone of $\Box$ at $\alpha$ and let $\sigma_\alpha$ be the
dual cone.

\begin{dfn}
\label{dfn:dual-edge-vectors}
I define the {\em dual edge vectors for $\Box$ at\/} $\alpha$ to be the
primitive generators of the cone $\sigma_\alpha$ in $N$. When $\sigma_\alpha$
is not basic, the dual edge
vectors do not generate $N$ as a group, but instead a
sublattice $N'_\alpha$ of $N$ of finite index, which I call the {\em dual edge
lattice for $\Box$ at\/} $\alpha$.
\end{dfn}

The cone $\sigma_\alpha$ is basic with respect to
$N'_\alpha$, and the corresponding variety $X'_\alpha=
X_{\sigma_\alpha,N'_\alpha}$ is smooth at
$P_\alpha$. By Corollary
\ref{prop:quotient}, the map
$X'_\alpha \to X_\alpha = X_{\sigma_\alpha,N}$ is the quotient map by the
action of the
finite abelian group
$K_\alpha=N/N'_\alpha\cong\hom_{\ZZ}(M'_\alpha/M,\QQ/\ZZ)$.  Here $M'_\alpha$
is the dual of $N'_\alpha$ and is naturally a superlattice of $M$. There is a
unique pairing $M'\times N \to \QQ/\ZZ$ which extends the pairings $M\times N
\to \ZZ$ and $M'\times N'\to \ZZ$. We then use the morphism $\QQ/\ZZ \to \CX$
given by the exponential map to identify $K_\alpha$ with
$\hom_{\ZZ}(M'_\alpha/M,\CX)$. If
we
identify $k\in K$ with the morphism $k:M'_\alpha\to\QQ/\ZZ$ such that
$k(M)=0$, the action is given by
\begin{equation}
\label{eq:Kaction}
k\cdot u'(m')=\exp(2\pi i\ip{k}{m'})u'(m'),
\end{equation}
for $u'\in \usig'$. Since the invariant part of $M'_\alpha$ under $K_\alpha$ is
$M$, the line bundles $L_\alpha$ and $L'_\alpha$ over $X_\alpha$ and
$X'_\alpha$ defined by the polytope $\Box$ are related by
$L_\alpha=L_{\alpha}'^K$. Equation (\ref{eq:coho}) shows that the cohomology of
$L_\alpha$ can be identified with the $K_\alpha$-invariant part of that of
$L_\alpha'$.

In summary, $(U'_\alpha, t, K_\alpha, L'_\alpha)$ is a local model for $X$ at
$P_\alpha$. Applying the Lefschetz formula for orbifolds, one deduces:

\begin{thm}
\label{thm:formula-sing}
For a simple convex polytope $\Box$ in $M$, we have
\begin{equation}
\label{eq:lef-fns}
\sum_{m\in\Box}m(t)=\sum_{\alpha\in \ext\Box}
\nu_\alpha(t)
\end{equation}
where
\begin{equation}\label{eq:nu-sing}
 \nu_\alpha(t)= {1\over{|K_{\alpha}|}}\sum_{k\in K_{\alpha}}
{\alpha(t)
\over
(1-e_k(\lambda_{\alpha}^{\prime 1})
\lambda_{\alpha}^{\prime 1}(t))
\cdots
(1-e_k(\lambda_{\alpha}^{\prime n})
\lambda_{\alpha}^{\prime n}(t))},
\end{equation}
and we write $e_k(\lambda)$ for $\exp(2\pi i{\ip{k}{\lambda}}).$ Here, the
vectors $\ntup{{\lambda_{\alpha}^{\prime
1}}}{\lambda_{\alpha}^{\prime n}}$ are the edge vectors of
$\Box$ at $\alpha$ in the dual $M'_\alpha$ of the dual edge lattice $N'_\alpha$
of definition \ref{dfn:dual-edge-vectors}, and $K_\alpha$ is the finite abelian
group $N/N'_\alpha$ acting according to equation (\ref{eq:Kaction}).
\end{thm}

\section{Laurent Expansions}
\label{sec:laurentexpansions}

In this section I expand the rational functions $\nu_\alpha$
away from their poles, i.e., in the domains where
$|\lambda_\alpha^i(t)|$ is not $1$, for $i=1,\dots,n$.
This has two benefits.

Firstly, it produces another formula which does not involve sums over roots of
unity. We shall use this in calculating the number of lattice points and the
volume.

Secondly it leads us to interpret the formula as a combinatorial statement,
decomposing the (characteristic polynomial for the) polytope $\Box$ as an
algebraic sum of the (characteristic series for the) cones $C_\alpha$ for each
extreme point. Ultimately this could be used to prove the formula using
elementary convex geometric reasoning. We don't attempt this here, as Ishida
has already reduced  the proof to the contractibility of convex sets
\cite{ishida}.

We begin by some general remarks about characteristic series for convex cones.

\subsection{Characteristic functions and series for convex cones}
 \label{subsec:characteristic}

We recall some notation, following \cite{ishida}. Let $A$ be a commutative ring
with identity. Recall that $A[M]$ denotes the {\em group algebra of $M$\/}
generated by elements $\bolde(m)$ for $m\in M$ satisfying relations
$\bolde(m)\bolde(m')=\bolde(m+m')$ and $\bolde(0)=1$. We denote by $A(M)$ the
total quotient ring of $A[M]$.

We define $A[[M]]={\rm Map}(M,A)$.   Elements $f\in A[[M]]$ can also be
expressed as formal Laurent series $f=\sum_{m\in M} f(m)\bolde(m)$ and this
defines a $A[M]$-module structure on $A[[M]]$ by:
$$\bolde(x)(\sum f(m)\bolde(m)) = \sum f(m-x) \bolde(m).$$

The relationship of $A[[M]]$ to $A(M)$ is as follows. To a given element
$\nu\in A(M)$ correspond (possibly) several elements of $A[[M]]$ called the
{\em Laurent expansions\/} of $\nu$. As we see below a convex cone $C$ in $M$
gives rise to elements $\nu^M_C\in A(M)$ and $\chi_{C\cap M}\in A[[M]]$ and the
latter is a Laurent expansion of the former.

\begin{dfn}
For $S$ a subset of $M$, we define the {\em characteristic series of $S$\/} to
be the element $\chi[S]=\chi_S$ of $A[[M]]$ corresponding to the set-theoretic
chacteristic function of $S$ (the function which takes values 1 on $S$ and 0
elsewhere), namely to the series $$\chi_S = \sum_{m \in S} \bolde(m).$$
\end{dfn}

Let $C$ be a (strongly convex rational simplicial) cone in $M_\RR$. We write
${\rm gen}^M_C=\{\lambda_1,\dots,\lambda_n\}$ for the primitive generators in
$M$ of $C$. The unit parallelepiped
$$Q^M_C=\{\sum a_i\lambda_i: 0\leq a_i < 1\}$$ defined by $C$ in $M$ intersects
$M$ in $\{c_1,\dots,c_{k}\}$. Here $k=|K|$, the order of the finite abelian
group which is the quotient of the dual lattice $N$ to $M$ by the lattice
generated by the primitive generators ${\rm
gen}^N_{C\dual}=\{\sigma^1,\dots,\sigma_n\}$ of $C\dual$ in $N$.

\begin{dfn}
For $C$ strictly convex, we define the {\em characteristic function for $C$
with respect to $M$\/}  is the following element of $A(M)$:
\begin{eqnarray*}
\nu^M_C & = & \sum_{c\in Q^M_C\cap M} \bolde(c) \prod_{\lambda\in{\rm gen}^M_C}
(1-\bolde(\lambda))\inv.\\
        & = & \sum_{j=1}^{|K|} \bolde(c_j) \prod_{i=1}^n
(1-\bolde(\lambda_i))\inv.
\end{eqnarray*}
For the translate of a cone $C$ by $\alpha\in M$, we set
$\nu^M_{\alpha+C}=\bolde(\alpha)\nu^M_C$.
\end{dfn}

Denote by ${\rm PL}_A(M)$ the $A[M]$-submodule of $A[[M]]$ generated by the set
of {\em polyhedral Laurent series\/}: $$\{\chi_{C\cap M} : C \mbox{ a basic
cone in }M_\RR\}.$$
 Ishida proves that the following \cite{ishida}
\begin{prop}
There exists a unique $A[M]$-homomorphism
$$\varphi:{\rm PL}_A(M) \to A(M)$$
such that $\varphi(\chi_{C\cap M})=\nu^M_C$, for all basic cones $C$ in
$M_\RR$.
\end{prop}

Actually,  we have:

\begin{prop}
For {\em any\/} cone $C$, $\chi_{C\cap M}\in {\rm PL}_A(M)$ and
$\varphi(\chi_{C\cap M})=\nu^M_C$ for $\varphi$ defined above.
\end{prop}
\begin{proof}
This follows from the remark that any element of $m\in M$ can be expressed
uniquely as $q+\sum x_i \lambda_i$ with $q\in Q^M_C\cap M$ and $x_i\in \NN$.
\end{proof}

The existence of $\varphi$ says essentially that we loose no information by
passing from the characteristic function of a cone to its' Laurent series, even
though the latter might not always have a well defined convergence on all of
$T_N$ (in the case $A=\CC$).

\paragraph{Remark} Whereas Ishida \cite{ishida} uses open cones, we find it
more convenient to use closed ones. The correpondence between the two is of
course that $C\cap M =\cup_{F < C} ({\rm int} F)\cap M$, where the union runs
over the faces of $C$.

\subsubsection{Action of $K$}

The group $K$ acts on $M'$ and hence on $A[[M']]$ by
$$k\cdot f = \sum_{m \in M} e_k(m)f(m)\bolde(m),$$ and we have $A[[M]] =
A[[M']]^K$. The following elementary remark gives the relationship between the
characteristic series for $C$ with respect to the two lattices $M$ and $M'$.

\begin{prop}
\label{prop:chi}
 For any cone $C$, we have
$$\chi_{C\cap M} = {1\over |K|} \sum_{k\in K} k\cdot \chi_{C\cap M'}.$$
\end{prop}
\begin{proof} Note that $k\cdot\chi(m')=e_k(m')\chi(m')$. Since $e_k$, for
$k\in K$, are nothing but the characters of the finite abelian group $M'/M$, we
have $e_k(M)=|K|$ and $e_k(m'+M)=0$, for all $m'\not\in M$. Hence the formula
follows.
\end{proof}

By the uniqueness of $\varphi$ we deduce that the same equality holds between
the characteristic functions of $C$:

\begin{cor} For any cone $C$, we have
$$\nu^M_C={1\over |K|} \sum_{k\in K} k\cdot \nu^{M'}_C.$$
\end{cor}

\subsection{Recovery of Brion's result}

We apply the results of the previous section with $A=\CC$. Then $\CC[M]$ is the
affine coordinate ring for the algebraic torus $T_N$ and its' field of
fractions $\CC(M)$ is the ring of rational functions on $T_N$.

The Lefschetz formula is expressing the chacteristic series $\chi_\Box$ of
$\Box$ as a sum of elements of $\CC(M)$. The theorem below says that these are
simply the characteristic functions for the tangent cones of $\Box$ at its'
extreme points. See \cite[Th\'eor\`eme 2.2]{brion}

\begin{thm}
\label{thm:formula-sing-b}
 Let $\Box$ be a simple convex polytope $\Box$ in $M$. Denote by $C_\alpha$
its' tangent cone at $\alpha\in\ext\Box$. Then we have \begin{equation}
\label{eq:lef-fns-b}
\chi_{\Box\cap M}= \sum_{\alpha\in\ext\Box} \nu^M_{C_\alpha}.
\end{equation}
\end{thm}
\begin{proof}
By theorem \ref{thm:formula-sing}  we have $\nu_\alpha= {1\over|K|}\sum_{k\in
K} k\cdot\nu^{M'}_{C_\alpha}$, which by the corollary of the previous section
is nothing but $\nu^M_{C_\alpha}$.
\end{proof}

\subsection{Laurent expansions of $\nu_C$ and their domains of validity}

We take $A=\CC$ and give all the different possible Laurent expansions of
$\nu^M_C$ for a cone $C$. When we attempt to evaluate these on elements of
$T_N$ these series only converge on certain open subsets which we specify here.

\subsubsection{The expansions}

We adopt the same notation as in section \ref{subsec:characteristic}. The
primitive generators of $C_\alpha$  are $\la^i$ in $M$ and
$\lambda_\alpha^{\prime i}$ in $M'_\alpha$.

\begin{prop}[Basic Expansion] For $|\lambda_\alpha^{\prime i}(t)|<1,$ for
$i=1,\dots,n,$ we have
\begin{equation}
\label{eq:nu-basic-expansion}
\nu_\alpha(t)=\chi_{\alpha+C_\alpha\cap M}(t).
\end{equation}
\end{prop}
\begin{proof}
Applying the elementary expansion (valid for $|z|<1$)
$$ (1-z)\inv= 1 + z + z^2 + z^3 + \cdots $$
to the individual factors $(1-e_k(\lambda_{\alpha}^{\prime
i})\lambda_{\alpha}^{\prime i}(t))\inv$ gives:
$$\nu_\alpha(t)=\alpha(t)
{1\over{|K_{\alpha}|}}\sum_{k\in K_{\alpha}}(\sum_{c_1,\dots,c_n=0}^{\infty}
 e_k(c\cdot \lambda'_\alpha)
(c\cdot \lambda'_\alpha)(t)),$$
where I have written $c\cdot \lambda'_\alpha$ for $\sum_{i=1}^{n}c_i
\lambda_{\alpha}^{\prime i}$. Since the series is convergent, one has
$$
\nu_\alpha(t) =
\sum_{c_1,\dots,c_n=0}^{\infty} (\alpha+c\cdot \lambda'_\alpha)(t)
{1\over{|K_{\alpha}|}}\sum_{k\in K_{\alpha}} e_k(c\cdot \lambda'_\alpha),$$
and the result follows from the proof of proposition \ref{prop:chi}.
\end{proof}

There are in fact $2^n$ different possible expansions for $\nu_\alpha(t)$
depending on whether we expand about $\lambda_\alpha^{\prime i}(t)=0$ or
$\infty$, each expansion being valid for
$ |\lambda_\alpha^{\prime i}(t)|<1$ or $>1$ respectively.

\paragraph{Notation:} Let $s$ be an $n$-tuple $s\in \{\pm1\}^n$. As a
shorthand, I will write:
\begin{eqnarray*}
\lambda'_\alpha  & \eqdef & (\lambda_\alpha^{\prime 1}, \dots,
\lambda_\alpha^{\prime n})\\
s\lambda'_\alpha & \eqdef & (s_1 \lambda_\alpha^{\prime 1}, \dots,
s_n\lambda_\alpha^{\prime n}).
\end{eqnarray*}

I also write $\langle\lambda'_\alpha\rangle$ for the cone
$\langle\lambda_\alpha^{\prime 1},\dots,\lambda_\alpha^{\prime n}\rangle$.
I define the quantity $s_{\minus}\cdot\lambda'_\alpha$ by:
$$s_{\minus}\cdot\lambda'_\alpha = \sum_{s_i=-1} s_i\lambda_\alpha^{\prime
i}.$$
An element $m\in M'$ defines a region $T_{m}$ of $T_{N'}$ by:
$$T_{m}=\{t\in T_{N'}: |m(t)|<1\}.$$
I also write, for a cone $C$ in $M$,
$$T_C=\{t\in T_{N'} : |m(t)|<1, \forall m\in C\cap M\}.$$
Thus, for example,
$$T_{\langle\lambda'_\alpha\rangle} = T_{\lambda_\alpha^{\prime
1}}\cap\cdots\cap T_{\lambda_\alpha^{\prime n}}.$$

\begin{prop}[General Expansion]
\label{prop:general-exp}
 Given  $s\in \{\pm1\}^n$, we have, for $t\in T_{\langle
s\lambda'_\alpha\rangle},$
\begin{equation}
\label{eq:nu-general-exp}
\nu_\alpha(t)=(\prod_{i=1}^n s_i)\chi[{\alpha + s_{\minus}\cdot\lambda'_\alpha
+ \langle s\lambda'_\alpha\rangle \cap M}](t).
\end{equation}
\end{prop}
\begin{proof}
In order to expand $\nu_\alpha$ when, for some $i$, we have
$|\lambda_\alpha^{\prime i}(t)|>1,$ I use the other expansion of $(1-z)\inv$,
valid for $|z|>1$:
$$(1-z)\inv= -z -z^2 -z^3 - z^4 - \cdots. $$
The result follows in the same way as the basic expansion. Note that compared
to the basic expansion, the cone whose characteristic series we end up with
undergoes  a reflection plus a translation: $\langle\lambda'_\alpha\rangle\cap
M$ becomes $s_{\minus}\cdot\lambda'_\alpha + \langle s\lambda'_\alpha\rangle
\cap M$. This is due to the shift from $1+z+z^2+\cdots$ to
$-z^1-z^2-z^3-\cdots$.
\end{proof}

\subsubsection{Consistency of expansions}

It doesn't make sense to expand all the $\nu_\alpha$ according to
(\ref{eq:nu-basic-expansion}) because the variable $t$ can't satisfy the
condition  $ |\lambda_\alpha^{\prime i}(t)|<1$ for all $i$ and $\alpha$. For
one thing, if $\alpha$ and $\beta$ are two extreme vertices of $\Box$ connected
by an edge, we will have
$\lambda_\alpha^{\prime i}=-\lambda_\beta^{\prime j}$ for some $i$ and $j$, so
that  $ |\lambda_\alpha^{\prime i}(t)|<1 \iff |\lambda_\beta^{\prime j}(t)|>1$.

I we can find a domain for $t\in T_{N'}$ such that {\em all\/} the expansions
we perform are valid {\em at the same time,} then when we sum up all the
$\nu_\alpha(t)$, all but a finite number of terms in the infinite series
cancel, and we get the characteristic polynomial $\chi_\Box$ evaluated on $t$.

For each $\beta\in\ext\Box$, we choose an element $s^\beta\in \{\pm1\}^n$, and
expand according to (\ref{eq:nu-general-exp}). We require that the set
\begin{equation}
\bigcap_{\beta\in\ext\Box} T_{\langle s^\beta \lambda'_\beta\rangle} =  T_{\cup
\{\langle s^\beta \lambda'_\beta\rangle : {\beta\in\ext\Box}\}}
\end{equation}
be non-empty. I turn next to the necessary conditions for this to be so.

\subsubsection{Neccessary conditions for a consistent expansion}
\label{subsub:necc-cond}

The above requirement implies, for instance, that if $\lambda^{\prime i}_\alpha
=-\lambda^{\prime j}_\beta$, as it happens for ajdacent vertices, then
$s^\alpha_i=-s^\beta_j$. This can be thought of graphically as choosing a
direction for each edge of the polytope $\Box$ and sticking to it throughout
the expansion. For each vertex $\alpha$ if the $i$-th edge is pointing into
$\alpha$ then we set $s^\alpha_i=-1$, if it is pointing out, we set
$s^\alpha_i=+1.$

Another necessary condition is that we choose $s^\alpha=(1,1,\dots,1)$ for some
$\alpha\in\ext\Box$. This can be seen easily, if one thinks for a moment of
decomposing $\chi_\Box$ as a sum of characteristic series for cones:
\begin{equation}
\label{eq:sum-chi}
\chi_\Box = \sum_{\beta\in\ext\Box} \pm \chi_{C'_\beta\cap M}
\end{equation}
where the cones $C'_\beta$ are obtained from the tangent cones $C_\beta$
eventually by the `reflection + translation' process prescribed in the general
expansion in proposition \ref{prop:general-exp} and the sign is determined by
the number of reflections specified by $s^\beta$. One of the cones involved
must be $C_\alpha$, for some $\alpha\in\ext\Box$. It will have all of its'
edges pointing outwards in the above orientation and will correspond to the
characteristic series $+\chi_{C_\alpha\cap M}$. I will call this the {\em base
vertex\/} for the expansion.

 The non-emptiness requirement above then implies that the following condition
on the orientations be satisfied:

\paragraph{Orientation condition} Let $\lambda^{\prime i}_\alpha$ for
$i=1,\dots p$ be any set of edges emanating from $\alpha$ that have been
oriented so that they are {\em all outgoing with respect to $\alpha$.\/ } Then
we require that for all $\beta\neq\alpha,$
 \begin{equation}
\label{eq:exp-condition}
\hbox{if }(\lambda^{\prime}_\beta)^j \in \pm \langle\lambda^{\prime
1}_\alpha,\dots,\lambda^{\prime p}_\alpha\rangle \hbox{ then }(s^\beta)_j=\pm
1.
\end{equation}

In words, this says that if an edge $\lambda^{\prime j}_\beta$ is a linear
combination, all of whose coefficients are of the same sign or zero, of
oriented edges $\lambda^{\prime i}_\alpha$ all going outwards from a given
vertex $\alpha$, then it should be oriented in the direction which includes it
in the cone spanned by these outgoing edges. This is because, if it were
oriented oppositely, it would mean that $T_{\langle\lambda^{\prime
1}_\alpha,\dots,\lambda^{\prime p}_\alpha\rangle} \cap T_{s^\beta_j
\lambda^{\prime j}_\beta} = \emptyset,$ since one cannot have both
$|\lambda^{\prime j}_\beta(t)|<1$ and $|-\lambda^{\prime j}_\beta(t)|<1.$

\subsubsection{Domain of validity of simultaneous expansions}

It is always possible to choose at least one orientation of the edges of $\Box$
which satisfies the orientation condition (\ref{eq:exp-condition}). Suppose we
have chosen such an orientation. For what values of $t\in T_N$ is it valid ? In
order to answer this, let us first make some remarks about the regions
$T_C\subset T_N,$ for $C$ a cone in $M$.

It is helpful, to describe $T_C$, to decompose $T_N$ as $CT_N\times H$,
corresponding to the Lie algebra decomposition $\bt_\CC=\bt\oplus
i\bt$. By identifying the second factor in the Lie algebra
decomposition with $N_\RR$, we have the exponential map
 $$N_\RR \stackrel{\exp}{\to} H.$$

\begin{lemma} If $C$ is a cone in $M$, then $T_C$ is given by
$$T_C=CT_N\times\exp(-{\rm int}(C\dual))\subset CT_N\times H.$$
\end{lemma}
\begin{proof}
The interior of $C\dual$ is the set of $n\in N_\RR$ such that $\ip{n}{c}>0,
\forall c\in C.$ Under the exponential map, the orbit $CT_N\times \{-n\}$
corresponds to an orbit of constant modulus strictly less than $1$.
\end{proof}

{}From this, we see that
$$\bigcap_{\beta\in\ext\Box} T_{\langle s^\beta \lambda'_\beta\rangle} =
CT_N\times\exp(-{\rm int}(\sigma)),$$
where
\begin{equation}
\label{eq:sigma}
\sigma=\left(\bigcup_{\beta\in\ext\Box} \langle s^\beta
\lambda'_\beta\rangle\right)^\vee.
\end{equation}
If we respect condition (\ref{eq:exp-condition}), we see that
$\bigcup_{\beta\in\ext\Box} \langle s^\beta \lambda'_\beta\rangle$ never
contains a whole subspace, so that $\sigma$ is non-zero. The expansion
determined by $s^\beta$ for $\beta\in \Box$ is thus valid in the region
$T_\sigma\subset T_N$ given by equation (\ref{eq:sigma}).

\subsection{Elementary convex geometric interpretation}

According to the work we have done in the previous sections, one can prove the
extreme point formula as follows:

Begin by orienting the edges of $\Box$ such as to respect condition
(\ref{eq:exp-condition}). This defines a cone (with a sign) for each extreme
vertex, according to proposition \ref{prop:general-exp}, and the algebraic sum
of their characteristic series should yield the characteristic polynomial for
the polytope $\Box$. If one can prove this for one admissible orientation of
the edges of $\Box$, then the formula for the characteristic functions follows
by the existence of Ishida's $\CC[M]$-homomorphism in the previous section.
This gives a proof of the formula involving only elementary convex geometry. We
won't bother with this, as Ishida \cite{ishida} already gives a proof which
reduces the problem to the contractibility of convex sets.

Instead we can deduce the following result in convex geometry:

\begin{thm}
\label{thm:chi-decomposition}
For all orientations $\{s^\alpha\}$ of the edges of $\Box$ satisfying the
orientation condition (\ref{eq:exp-condition}) we have
$$\chi_{\Box\cap M} = \sum_{\alpha\in\ext\Box} \pm\chi_{C^s_\alpha\cap M}$$
where $\pm=\prod_i (s^\alpha)_i$ and
$$ C^s_\alpha = \alpha + s_{\minus}\cdot\lambda'_\alpha + \langle
s\lambda'_\alpha\rangle.$$
\end{thm}

\section{Number of Lattice Points and Volume}

In this section I expand the functions $\nu_\alpha(t)$
around $t=1$ and derive formulae for the number of lattice
points and volume of $\Box$.

\subsection{The Number of Lattice Points}

Equation (\ref{eq:lef-fns-b}) expresses an equality between
the finite Laurent polynomial determined by $\Box$ and a sum
a rational functions. When evaluated on $t\in T_N$ with $t\to 1$ the left-hand
side tends to the number of lattice points of $\Box$ whereas on the right-hand
side the rational functions may have poles.

I choose a one-parameter subgroup $\{\exp(s\zeta) : s\in\RR\}$ determined by
some element $\zeta$ of the Lie algebra $\bf t$ of $CT_N$.
Substituting $\exp(s\zeta)$ for $t$, the formula reduces to an equality between
rational functions of $s$ --- provided I choose a one-parameter subgroup that
does not coincide with the singular loci of the $\nu_\alpha$.

\begin{dfn} For short, I call $\zeta\in \bt$ {\em generic\/} if
$\ip\zeta{\la^i}\ne0$, for all $i$ and $\alpha$. (This is indeed the case
generically).
\end{dfn}

 For generic $\zeta$, the functions $\nu_{\alpha,\zeta}^\Box: s \mapsto
\nu_\alpha^\Box(e^{s\zeta})$
can expanded in Laurent series:
$$ \nu_{\alpha,\zeta}^\Box(s)=\sum_{i=-\infty}^{\infty}
\nu_{\alpha,\zeta,i}^\Box s^i,
$$
 and their sum as $s\to 0$ is
obviously given by the sum of the constant terms
$\nu_{\alpha,\zeta,0}$ in each expansion.

Denote by $C_\alpha$ the tangent cone of $\Box$ at $\alpha\in\ext\Box$, and by
$\lambda^i_\alpha$ for $ i=1,\dots,n$, its' primitive generators in $M$. The
semi-open unit parallelepiped determied by the generators of $C_\alpha$ in $M$
is denoted
\begin{equation}
\label{eq:Qalpha}
Q_\alpha=Q^M_{C_\alpha}=\{\sum a_i\la^i: 0\leq a_i < 1\}.
\end{equation}
We have
$$\nu^\Box_{\alpha,\zeta}(s)= {\sum_{q\in Q_\alpha\cap M}
e^{s\ip{\zeta}{\alpha+q}}
\over
(1-e^{s\ip\zeta{\la^1}})\cdots(1-
e^{s\ip\zeta{\la^n}} )},$$
 provided $\ip\zeta{\la^i}\ne0$.
The zero-th order term in the expansion of $\nu^\Box_{\alpha,\zeta}(s)$ is a
homogeneous function of $\zeta$, which is equal to:
$$ {\sum_{q\in Q_\alpha\cap M} e^{s\ip{\zeta}{\alpha+q}}\over s^n\prod_i
(-\ip\zeta{\la^i})}{\prod_i
(-s\ip\zeta{\la^i})\over\prod_i(1-\exp(s\ip\zeta{\la^i}))}
$$
which gives
$${1\over\prod_i \ip\zeta{\la^i}}
\sum_{j=0}^n {(-1)^j \over j!} \sum_{q\in Q_\alpha\cap M}
\ip{\zeta}{\alpha+q}^j {\cal T}_{n-j}(\ip\zeta\la),
$$
 where ${\cal T}_k $ are the {\em Todd polynomials}, homogeneous polynomials of
degree $k$ whose coefficients can be expressed in terms of the Bernoulli
numbers \cite{hirz}. They are defined by the formal series
$$\sum_{k=0}^\infty s^k {\cal T}_k(x_1,x_2,\dots) = \prod_{i\geq 1}
{sx_i\over{1-\exp(-sx_i)}}.
$$
By ${\cal T}_k(\ip\zeta\la)$ I mean
$T_k(\ip\zeta{\la^1},\dots,\ip\zeta{\la^n})$.

\begin{thm}\label{thm:number}
 Let $\Box$ be a simple convex lattice polytope.
Denote by $C_\alpha$ the tangent cone of $\Box$ at $\alpha\in\ext\Box$, and by
$\lambda^i_\alpha,$ for $ i=1,\dots,n$, the primitive generators of $C_\alpha$
in $M$.
 The semi-open unit parallelepiped determied by the generators of $C_\alpha$ in
$M$ as in equation \ref{eq:Qalpha} is denoted $Q_\alpha$. Then, for generic
$\zeta\in \bt$, the number of lattice points in $\Box$ is given by
$$\sum_{\alpha\in\ext \Box}{1\over\prod_i \ip\zeta{\la^i}}
\sum_{j=0}^n {(-1)^j \over j!} \sum_{q_\alpha\in Q_\alpha\cap M}
\ip{\zeta}{\alpha+q_\alpha}^j {\cal T}_{n-j}(\ip\zeta\la).$$
\end{thm}

\paragraph{Remark 1} It might be more convenient in some cases to subdivide the
tangent cone into non-singular cones. One obtains a similar formula (see
\cite[Th\'eor\`eme 3.1]{brion}).

\paragraph{Remark 2} Putting
$t=\exp(s\zeta)$ corresponds to considering the Lefschetz number for the action
of the one-parameter subgroup $G_\zeta$ of $CT_N$ generated by $\zeta\in\bt =
{\rm Lie }\,CT_N$. Generically this has a dense orbit, and therefore the same
fixed points on $X$ as the whole real torus $CT_N$, and so the Lefschetz
formula for $G_\zeta$ is the same as that obtained by substituting
$\exp(s\zeta)$ for $t$.
This is not true of course when $\ip\zeta{\la^i}=0$, for some $i$ and $\alpha$.
Indeed in that case the group $G_\zeta$ has whole circles of fixed points.
Restricting to $G_\zeta$ corresponds to projecting the vertices and edges of
$\Box$ onto the hyperplane in $M_\RR$ defined by the form $\zeta\in
N_\RR$.

\subsection{The Volume}

\subsubsection{The ``Classical Limit''}

In the introduction I mentioned the fact that for larger and larger
polytopes (or finer and finer lattices) the number of points is
asymptotically equal to their volume --- I call this ``the classical limit'' by
analogy with the limit $\hbar\to 0$ in quantum mechanics. More precisely, for
any $n$-dimensional polytope $\Box$, the volume of $\Box$ is given by
\begin{equation}
\label{eq:vol_lim}
{\rm vol}_n(\Box)=\lim_{k\to\infty}{\#(k\inv M\cap
\Box)\over k^n} = \lim_{k\to\infty}{\#(M\cap k\Box)\over k^n}.
\end{equation}
 Indeed \cite{mac:poly}, the function
$$H_{\Box}(k)=\#(k\inv M\cap \Box)=\#(M\cap (k\Box))$$
is a polynomial of degree $n$, for $k\in\NN$, with leading
coefficient ${\rm vol}_n(\Box),$ and is called the {\em
Hilbert polynomial\/} for $\Box$. The polynomial $H_{\Box}$ is in fact equal to
the {\em Hilbert polynomial\/} $H_{(\XB,\LB)}$ for the pair $(\XB,\LB)$, namely
$$H_{(\XB,\LB)}(k)=\chi(\XB,{\cal O}_{\XB}(k\LB))=
\sum(-1)^i\dim H^i(\XB,{\cal O}_{\XB}(k\LB)).$$
 This follows from equation  (\ref{eq:coho}) and because taking tensor powers
$\LB^{\otimes k}$ of $\LB$ corresponds to taking multiples of $kN$ of $N$, and
hence submultiples $k\inv M$ of $M$.

 \begin{thm}
\label{thm:volume} Let $\Box$ be a simple convex lattice polytope and adopt the
same notation as theorem \ref{thm:number}. Let $|K_\alpha|$ denote the order of
the singularity of $\Box$ at $\alpha$. Then for generic $\zeta$ the volume of
$\Box$ is given by
$${\rm vol}_n(\Box)= {(-1)^n\over n!}
\sum_{\alpha\in\ext\Box}{\ip\zeta{\alpha}^n |K_\alpha| \over
\ip\zeta{\la^1}\cdots\ip\zeta{\la^n}}.$$
\end{thm}
\begin{proof}  The proposition follows from taking the coefficients
 of the $k^n$ terms in theorem \ref{thm:number} applied to the
 polytope $k\Box$. Note that $\ext k\Box=k(\ext\Box)$ and that
$C^{k\Box}_{k\alpha}=C^\Box_\alpha.$ Note that the order $|K_\alpha|$ of the
singularity at $\alpha$ is equal to the cardinality of $Q_\alpha\cap M$. See
\cite[Corollaire 2]{brion}.
\end{proof}

\subsubsection{The Riemann-Roch approach}

The volume of $\Box$ appears if one uses the same geometric approach based on
the $d''$-complex but directly applies the Riemann-Roch theorem, instead of
computing the Lefschetz number for the action of $t\in T$ and then letting
$t\to1$.

The Riemann-Roch theorem expresses the Euler characteristic of a
holomorphic vector bundle $E$ over a complex manifold $X$ in terms of
characteristic classes of $E$ and (tangent bundle to the) $X$:
\begin{equation}
\chi(X,E)=\{\hbox{ch}(E)\cdot{\cal T}(X)\}[X],
\end{equation}
 where ch$(E)$ and ${\cal T}(X)$ are the Chern character of $E$
and the Todd class of $X$, respectively. If $E$ has rank $n$ and
$c_1,\dots,c_n$ denote the characteristic classes of $E$ then the {\em Chern
character} can be defined by the power series
$$\sum_{i=1}^n e^{x_i}= n+\sum x_i+{\sum x_i^2\over 2!}+\cdots,$$
where the $c_i$ are to be thought of formally as the elementary symmetric
functions in the $x_i$.

Since we are in a one-dimensional situation and $c_1(\LB)$ is represented by
the K\"ahler form $\omega$, the Chern character is given by
$${\rm
ch}(\LB)=1+\omega+{\omega^2\over2!}+{\omega^3\over3!}+\dots+{\omega^n\over
n!}.$$
The {\em Todd class} is a polynomial in the characteristic
classes $c'_i$ of the tangent bundle of $X$. If the $c'_i$ are regarded
formally as the elementary symmetric functions of the $x'_i$ (as in the case
above), the Todd class can be expressed as
 $${\cal T}(X)=\prod_i {x'_i\over 1-e^{-x'_i}}.$$
(Presumably, there is some relationship between these and the Todd polynomials
of theorem \ref{thm:number} which in this case     exhibits the Riemann-Roch
formula as the ``classical limit'' of the Lefschetz fixed point formula.)
By multiplying the two series selecting the terms of order $n$ and evaluating
them on
$[X]$, we get
$$\chi(X,\LB) ={\rm vol}_n(X)+ \hbox{\em lower order terms},$$
where the ``lower order terms" are terms involving powers of $\omega$ of order
less than $n$. Again, because refining the lattice $M$ corresponds to
multipying $\omega$, we see that $\chi(X,t\LB)$ is given asymptotically by
${\rm vol}_n(X)t^n$.

%Bibliography for Toric Varieties

\end{document}